\documentclass[reprint,amsmath,author-numerical]{revtex4-1}
\usepackage{amssymb,amsthm,tikz,pgfplots}
\usepgfplotslibrary{groupplots}

\usepackage{graphicx}
\graphicspath{}
\usepackage{float}
\usepackage{dcolumn}
\usepackage{bm}
\usepackage{subfig}
\usepackage{enumerate}
\usepackage{color,transparent} 
\usepackage{bbm}
\usepackage{multirow}
\usepackage{hyperref}
\providecommand{\eps}{\varepsilon}
\providecommand{\R}{\mathbb{R}}
\providecommand{\N}{\mathbb{N}}
\providecommand{\Z}{\mathbb{Z}}
\providecommand{\cO}{\mathcal{O}}
\providecommand{\cA}{\mathcal{A}}
\providecommand{\cL}{\mathcal{L}}
\providecommand{\cP}{\mathcal{P}}
\providecommand{\argmin}{{\rm argmin}}
\providecommand{\spn}{{\rm span}}

\providecommand\intr[1]{\mathring{#1}}
\providecommand{\nin}{\ell}
\providecommand{\nbd}{m}
\providecommand{\inn}{\mathbf{in}}
\providecommand{\bd}{\mathbf{bd}}

\providecommand{\OJ}[1]{{#1}}

\usepackage{listings}
\definecolor{colKeys}{rgb}{0,0,1}
\definecolor{colIdentifier}{rgb}{0,0,0}
\definecolor{colComments}{rgb}{0,1,0.3}
\definecolor{colString}{rgb}{0,0.5,0}

\definecolor{dkgreen}{rgb}{0,0.6,0}
\definecolor{gray}{rgb}{0.5,0.5,0.5}
\definecolor{lightgray}{rgb}{0.9,0.9,0.9}

\newfloat{Code}{H}{myc}

\providecommand{\evplot}[4]{
\begin{tikzpicture}
    \begin{axis}[enlargelimits=false,
    	axis on top,
		width=0.25\textwidth,
		height=0.25\textwidth,
		xtick={0,6.28},ytick={0,6.28},
		xticklabels={0,$2\pi$},yticklabels={\empty,$2\pi$},
		title={#2},
		xlabel={#3},xlabel style={at={(0.5,0.1)}},
		ylabel={#4},ylabel style={at={(0.3,0.5)}},
		title style={font=\footnotesize}]
        \addplot graphics [xmin=0,xmax=6.28,ymin=0,ymax=6.28]
			{#1};
    \end{axis}
\end{tikzpicture}}

\providecommand{\cylevplot}[4]{
\begin{tikzpicture}
    \begin{axis}[enlargelimits=false,
    	axis on top, axis equal image,
		width=0.28\textwidth,
		height=0.28\textwidth,
		xtick={0,6.28},ytick={0,3.14},
		xticklabels={0,$2\pi$},yticklabels={\empty,$\pi$},
		title={#2},
		xlabel={#3},xlabel style={at={(0.5,0.2)}},
		ylabel={#4},ylabel style={at={(0.3,0.5)}},
		title style={font=\footnotesize}]
        \addplot graphics [xmin=0,xmax=6.28,ymin=0,ymax=3.14]
			{#1};
    \end{axis}
\end{tikzpicture}}

\begin{document}

\title{On fast computation of finite-time coherent sets using radial basis functions}

\author{Gary Froyland}
\affiliation{
School of Mathematics and Statistics, The University of
New South Wales, Sydney NSW 2052, Australia
}%
\thanks{The research of GF is supported by an ARC Future Fellowship (FT120100025) and an ARC Discovery Project (DP150100017).  GF and OJ also thank the Mathematisches Forschungsinstitut Oberwolfach for providing excellent working conditions for the initial stages of this research.}
\author{Oliver Junge}
\affiliation{
Zentrum Mathematik - M3
Technische Universit\"{a}t M\"{u}nchen
85747 Garching bei M\"{u}nchen, Germany
}

\date{\today}

\begin{abstract}
Finite-time coherent sets inhibit mixing over finite times.  The most expensive part of the transfer operator approach to detecting coherent sets is the  construction of the operator \OJ{itself}.
We present a numerical method based on radial basis function collocation and apply it to a recent transfer operator construction \cite{F15} that has been designed specifically for purely advective dynamics.
The construction \cite{F15} is based on a ``dynamic'' Laplace operator and minimises the boundary size of the coherent sets relative to their volume.
The main advantage of our new approach is a substantial reduction in the number of Lagrangian trajectories that need to be \OJ{computed}, leading to large speedups in the transfer operator analysis \OJ{when this computation is costly}.
\end{abstract}


\maketitle

\section{Introduction}
Finite-time coherent sets (FTCS) \cite{FSM10,F13,FPG14} in nonlinear, possibly time-dependent, dynamical systems are connected regions that are maximally dynamically disconnected from surrounding phase space when evolved over a specified time duration of finite length.
If the dynamical system is an advection-diffusion equation, e.g.\ a Fokker-Planck equation, finite-time coherent sets are those regions for which there is minimal exchange of a passive tracer across their boundary.
On the other hand if the dynamical system is a model of (purely advective) fluid flow, finite-time coherent sets are those regions for which there is minimal exchange of fluid across their boundary in the presence of small noise or diffusion.
The exchange of fluid across a boundary by small diffusion is proportional to the size of the boundary, and in the case of purely advective fluid flow, one can devise an alternative characterization of FTCS which minimizes boundary size relative to volume \cite{F15}.

In the pure advection setting, for volume-preserving dynamics, the constructions in this alternative characterization have been shown to arise as zero-diffusion limits of the ``classical'' FTCS approach \cite{F13}.
Thus, the approach \cite{F15} can be viewed both as an equivalent
geometric reinterpretation of the classical mathematical definitions based on probability and diffusion \cite{F13}, and as an alternative numerical method for identifying FTCS.
In the present paper we apply numerical analysis methods based on radial basis function collocation to very efficiently compute estimates of the main transfer operator object \cite{F15}.
Because of the flexibility of radial basis functions, this approach is suited to irregularly-shaped domains and computations based on trajectory data.

Denote by $M$ our domain of interest at the initial time $t_0$.
The domain $M$ is a compact, connected subset of $\mathbb{R}^d$ with $C^1$ boundary (formally, for the theorems quoted later, we will assume $M$ is in fact a compact, connected, $C^\infty$ Riemannian manifold of vanishing curvature).
Our advective dynamics is either a nonlinear volume-preserving map $T:M\to T(M)$ (or several compositions of possibly different maps), or a nonlinear volume-preserving flow map $\Phi(x,t_0,t)$, which solves the ODE $\dot{x}=F(x,t)$.
In the latter case, we assume that $F(x,t)$ is defined on a large enough domain to advect our initial domain $M$ forward for some finite time.

In the purely advective setting, how well a set ``mixes'' is generally related to the length or irregularity of the boundary of the set.
Measures of mixing specifically designed for pure advection, such as the mix-norm \cite{MathewMezicPetzold,DoeringThiffeault,Thiffeault} penalize the geometric irregularity of advected passive scalar concentrations by comparing the difference of the concentration from a uniform concentration over an infinite range of spatial scales.
Other advective mixing measures directly measure the size of the boundary of an initially compactly-supported passive scalar concentration \cite{KrotterChristovOttinoLueptow}.

\begin{figure}[htbp]
  \centering
  \includegraphics[width=0.4\textwidth]{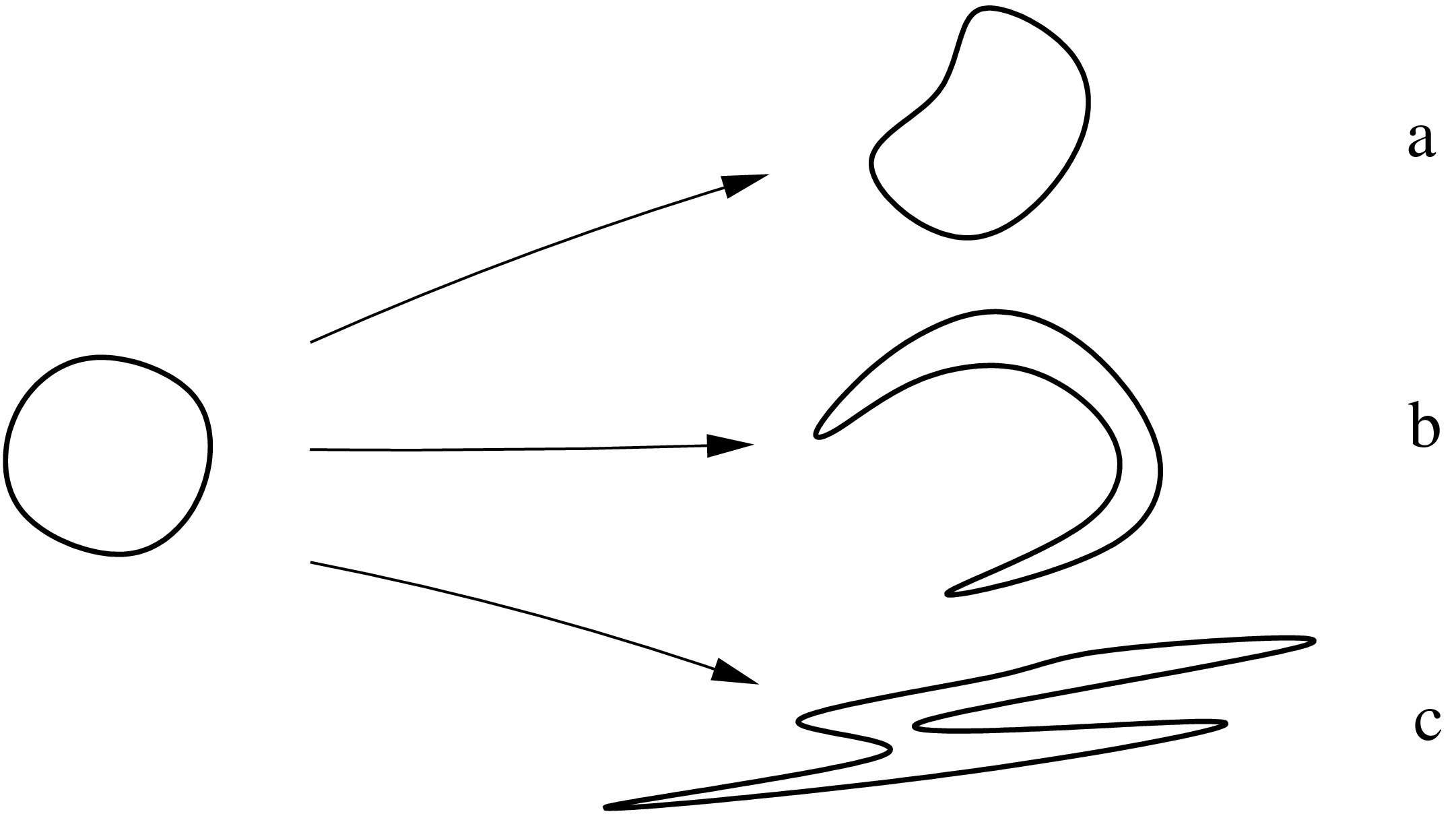}\\
  \vspace{.5cm}
  \caption{The two-dimensional shape on the left has a low boundary size to volume ratio.  The shapes on the right are three possible images of the shape on the left under three different nonlinear volume-preserving dynamical systems over a finite time duration.  Under dynamics `a' the set on the left retains a low boundary size to volume ratio, but under dynamics `b' and `c' the boundary size is significantly increased.}\label{fig:coherence}
\end{figure}

Our objective is to identify regions of $M$ that have low boundary size relative to volume, \emph{and retain} low boundary size relative to volume \emph{as they are advected}.
Such regions do \emph{not} evolve to have highly filamented boundaries, which under the addition of small noise would lead to substantial mixing; see Figure \ref{fig:coherence}.
By minimising the size of the boundary relative to volume either throughout the flow duration or at the beginning and the end of the flow, we find finite-time coherent sets, which mix least with the surrounding phase space.
\vspace{.1cm}

In recent years there have been several geometric methods proposed to characterize either trajectories or co-dimension 1 surfaces that represent coherent structures \cite{mancho_M,mezicmeso,Haller_11,Haller_12,allshouse_thiffeault,ma_bollt_shape,romkedar_M,FPG15}.
In some two-dimensional cases \cite{Haller_11,Haller_12}, the aim of the geometric methods is to find curves that undergo small amounts of stretching.
The new approach based on the dynamic Laplacian \cite{F15} directly minimizes boundary size (in arbitrary dimensions) when subjected to finite time dynamics.
Our main contribution in the present work is to reduce the computational effort by taking advantage of the underlying functional representation of the transfer operator approach.

Our approach is similar in spirit to the work of Williams {et al.}\cite{WiRyRo14a} in this issue.  There, the authors use a Galerkin scheme with globally supported radial basis functions to compute the constructions in \cite{F13}.
Here, we use collocation, locally supported functions, and the diffusion-free operator proposed by Froyland \cite{F15} (rather than the diffusive approach \cite{F13}) -- but the goal is very much the same, namely constructing (data-)efficient discretisations of the relevant operators in order to compute coherent sets.
We believe that ultimately, locally supported functions might yield the more efficient numerical scheme when aiming for high resolution.

\section{Advective setup for finite-time coherent sets}

For $d\ge 1$ we denote $d-1$-dimensional volume by $\ell_{d-1}$;  for example, if our fluid is three-dimensional $\ell_2$ measures surface area, while if our fluid is two-dimensional, $\ell_1$ measures curve length.
All of the theory and techniques outlined here work in arbitrary finite dimensions $d\ge 1$.

Let $\Gamma\subset M$ be a smooth $d-1$-dimensional object that separates $M$ into two subregions $M_1$ and $M_2$.
The hypersurface $\Gamma$ is the common boundary of $M_1$ and $M_2$, whose size we wish to minimise relative to the smaller of the volumes $\ell_d(M_1),\ell_d(M_2)$.
Denote by $\Phi(\cdot,t_0,t)$ the flow map from time $t_0$ to time $t$.  Below we define an objective function that calculates the combined boundary size to volume ratios at an initial time $t_0$ and a final time~$t_\mathbf{f}$.

\vspace*{0.2in}
\textbf{Definition 1:}
\label{cheegerdefn}
Define
\begin{equation}
\label{cheegereqn}
\mathbf{h}^D_{\{t_0,t_\mathbf{f}\}}(\Gamma)=\frac{1}{2}\left(\frac{\ell_{d-1}(\Gamma)+\ell_{d-1}(\Phi(\Gamma,t_0,t_\mathbf{f}))}{\min\{\ell_d(M_1).\ell_d(M_2)\}}\right),
\end{equation}

We wish to find the minimising $\Gamma$:  $\Gamma^*=\argmin_\Gamma \mathbf{h}^D_{\{t_0,t_\mathbf{f}\}}(\Gamma)$.
The finite-time coherent set is defined to be whichever of $M_1$, $M_2$ has the least volume;  the boundary is formed by the hypersurface $\Gamma^*$.

The above mathematical construction in the absence of any dynamics arises naturally in \emph{isoperimetric theory}\cite{chavelisoperimetric}, where purely geometric questions of how to disconnect a compact manifold in such a way that the disconnecting hypersurface size relative to the volume of the smallest disconnected piece is minimised.
Answers to such classical questions reveal important information about the geometry of the manifold.
On the other hand, in the present constructions the nonlinear dynamics plays the dominant role.
We refer the reader to \cite{F15} for further further background and discussion of these ideas.

If the flow time is long, it could be that an evolved region develops a filamented boundary during some intermediate time interval, but then loses this filamentation by the final flow time.
For the times that the filament is present, small-scale diffusion will decrease the coherence of the set.
In some applications, it may be important to not allow such transient filaments.
By extending the above basic characterisation of FTCS to include a series of times $t_0,t_1,\ldots,t_\mathbf{f}$, one can penalise such transient filamentations.
Define
\begin{equation}
\label{hdynn}
\mathbf{h}_{\{t_0,t_1,\ldots,t_{n-2},t_\mathbf{f}\}}^{D}(\Gamma)=\frac{\frac{1}{n}\sum_{i=0}^{n-1}\ell_{d-1}(\Phi(\Gamma,t_0,t_i))}{\min\{\ell_d(M_1),\ell_d(M_2)\}},
\end{equation}
where $t_{n-1}:=t_\mathbf{f}$,
as the natural generalisation of $\mathbf{h}^D_{\{t_0,t_\mathbf{f}\}}(\Gamma)$.
Alternatively, for continuous penalisation, we can consider
\begin{equation}
\label{hdynt}
\mathbf{h}_{[t_0,t_\mathbf{f}]}^{D}(\Gamma)=\frac{\frac{1}{\tau}\int_0^\tau \ell_{d-1}(\Phi(\Gamma,t_0,t))\ dt}{\min\{\ell_d(M_1),\ell_d(M_2)\}},
\end{equation}
as a time-continuous generalisation of $\mathbf{h}^D_{\{t_0,t_\mathbf{f}\}}(\Gamma)$;  see \cite{F15} for further details.

\section{A dynamic Laplace operator}

We use eigenfunctions of a dynamic Laplace operator to numerically find the minimising hypersurface $\Gamma^*$.
To motivate this approach, we introduce the \emph{dynamic Sobolev constant} \cite{F15}.
Define
\begin{equation}
\label{soboleveqn}
\mathbf{s}^D_{\{t_0,t_\mathbf{f}\}}(f)=\frac{\|\nabla f\|_{L^1}+\|\nabla(f\circ \Phi(\cdot,t_\mathbf{f},t_0))\|_{L^1}}{2\inf_{\alpha\in \mathbb{R}} \|f-\alpha\|_{L^1}}.
\end{equation}
The above constant is modelled on the Sobolev constant common in isoperimetric theory \cite{chavelisoperimetric}, where only the first term in the numerator is present;  see \cite{F15} for details.
There is a very strong formal connection between the constants $\inf_\Gamma\mathbf{h}^D_{\{t_0,t_\mathbf{f}\}}(\Gamma)$ (purely geometric) and $\inf_{f\in C^\infty(M)}\mathbf{s}^D_{\{t_0,t_\mathbf{f}\}}(f)$ (purely functional), namely
\begin{equation}
\label{dffeqn}
\inf_\Gamma\mathbf{h}^D_{\{t_0,t_\mathbf{f}\}}(\Gamma)=\inf_{f\in C^\infty(M)}\mathbf{s}^D_{\{t_0,t_\mathbf{f}\}}(f)
 \end{equation}
 (see \cite{F15} for the formal statement).
This \emph{dynamic Federer-Fleming Theorem} is a generalisation of the celebrated Federer-Fleming Theorem \cite{federerfleming} known throughout differential geometry.

Our numerical approach will be to find an $f$ minimising $\mathbf{s}^D_{\{t_0,t_\mathbf{f}\}}(f)$.
In general, the infimum in the RHS of (\ref{dffeqn}) is only approached by a sequence of $C^\infty$ functions, and there is no simple formula for such a sequence of functions.
On the other hand if we use the $L^2$ norm instead of the $L^1$ norm in the RHS of (\ref{dffeqn}), the infimum is attained by a smooth function, which
is the second eigenfunction of the dynamic Laplace operator, defined by
\begin{equation}
\label{hattriangle}
\triangle^D:=(\triangle+\mathcal{P}^*\triangle\mathcal{P})/2.
\end{equation}
In (\ref{hattriangle}), $\triangle$ is the standard Laplace operator, $\mathcal{P}f=f\circ \Phi(\cdot,t_\mathbf{f},t_0)$ is the Perron-Frobenius (or transfer) operator for $\Phi(\cdot,t_0,t_\mathbf{f})$, and $\mathcal{P}^*f=f\circ \Phi(\cdot,t_0,t_\mathbf{f})$ is its dual, the Koopman operator for $\Phi(\cdot,t_0,t_\mathbf{f})$.
Note that the first Laplace operator in (\ref{hattriangle}) acts on $M$, while the second Laplace operator acts on $\Phi(M,t_0,t_\mathbf{f})$;  this is crucial for $\triangle^D$ to feel the geometry at both time instants.
In fact, $\mathcal{P}^*\triangle\mathcal{P}$ is the natural pullback of the Laplace operator on $\Phi(M,t_0,t_\mathbf{f})$ under $\Phi(\cdot,t_0,t_\mathbf{f})$, and is the Laplace-Beltrami operator on $M$ with respect to the pullback of the Euclidean metric on $\Phi(M,t_0,t_\mathbf{f})$.
Let $\delta$ denote the trivial (Euclidean) Riemannian metric on $\Phi(M,t_0,t_\mathbf{f})$; pulling $\delta$ back under $\Phi(\cdot,t_0,t_\mathbf{f})$ we obtain the Riemannian metric $\Phi(\cdot,t_0,t_\mathbf{f})^*\delta$ on $M$, and the map $\Phi(\cdot,t_0,t_\mathbf{f}):(M,\Phi(\cdot,t_0,t_\mathbf{f})^*\delta)\to (\Phi(M,t_0,t_\mathbf{f})),\delta)$ is an isometry.
Then $\triangle_{\Phi(\cdot,t_0,t_\mathbf{f})^*\delta}f=(\triangle_\delta(f\circ \Phi(\cdot,t_0,t_\mathbf{f}))^{-1})\circ \Phi(\cdot,t_0,t_\mathbf{f})=\mathcal{P}^*\triangle_\delta\mathcal{P}$, where $\triangle_\delta$ is the Laplace-Beltrami operator on the Riemannian manifold $(\Phi(M,t_0,t_\mathbf{f})),\delta)$ see Section 3.2 \cite{F15} or p27 \cite{chaveleigenvalues}.

If $M$ has a boundary, then the eigenfunction $f:M\to\mathbb{R}$ is required to satisfy generalised (oblique) Neumann boundary conditions.
Denote by $\mathbf{n}(x)$ the unit outward normal at $x\in\partial M$.
We require
\begin{equation}
\label{strongbc0}
\nabla f(x)\cdot\left[(I+C_{x,t_0,t_\mathbf{f}}^{-1})\mathbf{n}(x)\right]=0,\qquad x\in\partial M,
\end{equation}
where $C_{x,t_0,t_\mathbf{f}}:=D\Phi(x,t_0,t_\mathbf{f})^\top\cdot D\Phi(x,t_0,t_\mathbf{f})$ is  the (right) Cauchy-Green deformation tensor for the transformation $\Phi(\cdot,t_0,t_\mathbf{f})$.
This boundary condition also has a natural pullback interpretation, see Section 3.2 \cite{F13}.

\subsection{Spectral properties of the dynamical Laplacian}
The spectral properties of the $L^2$ eigenproblem (\ref{hattriangle})-(\ref{strongbc0}) are (see \cite{F15} for formal statements):
\begin{enumerate}
\item $\triangle^D$ is self-adjoint.
\item The eigenvalues form a decreasing sequence $0=\lambda_1\ge \lambda_2\ge\cdots$ with $\lambda_n\to-\infty$.
\item The leading eigenfunction $f_1\equiv 1$, and the eigenfunctions $f_1,f_2,\ldots$ corresponding to distinct eigenvalues are pairwise orthogonal in $L^2$.
\item One has the following variational characterisation of eigenvalues: if $f_1,f_2,\ldots$ are arranged to be orthonormal, denoting $X_k=\spn\{f_1,f_2,\ldots,f_k\}$
\begin{eqnarray}
\nonumber\lambda_k &=&
-\inf_{f}\left\{\frac{\|\nabla f\|_{L^2}^2+\|\nabla(f\circ\Phi(\cdot,t_\mathbf{f},t_0))\|^2_{L^2}}{2\|f\|_{L^2}^2}\right.\\
\label{variational}
&&\left.\qquad:\langle f,f_i\rangle=0, i=1,\ldots,k-1.\right\}
\end{eqnarray}
with the infimum achieved only when $f=f_k$.
Note that setting $k=2$ minimises (\ref{soboleveqn}) when $\|\cdot\|_{L^1}$ is replaced with $\|\cdot\|_{L^2}$.
\end{enumerate}

From (\ref{variational}), one sees eigenvalues $\lambda_k$ that are far from zero correspond to functions $f$ that have high gradient over a large part of the domain (we call such functions ``irregular'' as they are far from the most regular function, namely the constant function).
The equality (\ref{dffeqn}) makes this exact for $k=2$, and in fact, one also has a dynamic Cheeger inequality\cite{F15}:
\begin{equation}
\label{dceeqn}
\inf_\Gamma\mathbf{h}^D_{\{t_0,t_\mathbf{f}\}}(\Gamma)\le 2\sqrt{-\lambda_2},
\end{equation}
which links the least possible dynamic boundary size to volume ratio with the first nontrivial eigenvalue $\lambda_2$ of the dynamic Laplace operator $\triangle ^D$.

If there are $K>1$ independent ways to disconnect $M$, such that the disconnections each have similarly small dynamic boundary size to volume ratio, then there should be a cluster of eigenvalues $\lambda_2,\ldots,\lambda_{K+1}$ that are much closer to 0 than $\lambda_{K+2},\lambda_{K+3},\ldots$.
In such a situation, the $K$ finite-time coherent sets can be extracted from level sets of the eigenfunctions $f_k$, $k=2,\ldots,K+1$.
Perhaps the simplest way to do this is to
find the $K$ level sets $\Gamma$ by optimising $\mathbf{h}^D(\Gamma)$ with $K$ separate line searches.
For related approaches in the autonomous dynamics setting, see also \cite{froyland-dellnitz_03,deuflhard-weber_05}.

\subsection{Coherent sets and objectivity}
Existing transfer operator methods for identifying finite-time coherent sets \cite{FSM10,FHRSS12,F13,FPG14} extract the coherent sets from super- or sub-level sets of singular vectors of the transfer operator corresponding to singular values close to 1.
We follow the same strategy here, using level sets of the solutions of the eigenproblem (\ref{hattriangle})--(\ref{strongbc0}) that correspond to eigenvalues close to 0.
The numerical algorithm is detailed in Section \ref{sect:extract}.

Concerning objectivity of the algorithm, it is shown in \cite{F15} that if the domain $M$ is additionally subjected to a time-dependent affine transformation with orthogonal linear part, the solutions of the corresponding eigenproblem (\ref{hattriangle})--(\ref{strongbc0}) are simply transformed versions of the solutions of the original eigenproblem.
Thus, as the coherent sets are extracted from level sets of the eigenfunctions of (\ref{hattriangle})--(\ref{strongbc0}), the extracted coherent sets are also identically transformed as required for objectivity of the method.

\subsection{Penalising boundary size at multiple times}

If one wishes to penalise the boundary size at several discrete times between $t_0$ and $t_\mathbf{f}$ then define
\begin{equation}
\label{bartrianglen}
\triangle^D_{t_0,t_1,\ldots,t_{n-2},t_\mathbf{f}}:=\frac{1}{n}\sum_{i=0}^{n-1}\mathcal{P}_{t_0,t_i}^*\triangle\mathcal{P}_{t_0,t_i},
\end{equation}
where $\mathcal{P}_{t_0,t_i}f=f\circ \Phi(\cdot,t_i,t_0)$.
Consider the eigenproblem
\begin{equation}
\label{strongeqnmultn}
\triangle^D_{t_0,t_1,\ldots,t_{n-2},t_\mathbf{f}}f(x)=\lambda f(x),\quad x\in \mathring{M},
\end{equation}
with boundary condition
\begin{equation}
\label{strongbcmultn}
\nabla f(x)\cdot\left[\sum_{i=0}^{n-1}C_{x,t_0,t_i}^{-1} \mathbf{n}(x)\right]=0,\quad x\in\partial M.
\end{equation}

If one wishes to penalise the boundary size at all times between $t_0$ and $t_\mathbf{f}$ then define
\begin{equation}
\label{bartrianglet}
\triangle^D_{[t_0,t_\mathbf{f}]}:=\frac{1}{\tau}\int_0^\tau\mathcal{P}_{t_0,t}^*\triangle\mathcal{P}_{t_0,t}\ dt,
\end{equation}
where $\mathcal{P}_{t_0,t}f=f\circ \Phi(\cdot,t,t_0)$.
Consider the eigenproblem
\begin{equation}
\label{strongeqnmultt}
\triangle^D_{[t_0,t_\mathbf{f}]}f(x)=\lambda f(x),\quad x\in \mathring{M},
\end{equation}
with boundary condition
\begin{equation}
\label{strongbcmultt}
\nabla f(x)\cdot\left[\int_0^\tau C_{x,t_0,t}^{-1} \mathbf{n}(x)\ dt\right]=0,\quad x\in\partial M.
\end{equation}

In both of the above cases, the obvious versions of (\ref{dffeqn}) and (\ref{dceeqn}) hold \cite{F15}.

\section{Discretization of the dynamic Laplacian eigenproblem}
\label{sect:RBF}

Following Platte and Driscoll \cite{PlDr03a}, we are going to approximate the eigenproblem (\ref{hattriangle})-(\ref{strongbc0}) by collocation with radial basis functions.  Here, we choose the \emph{Wendland} functions  $\psi=\psi_{d,k}:[0,\infty)\to\R$ for various $k\in\N$, which have support on $[0,1]$, are polynomials of a certain degree and lead to strictly positive definite interpolation matrices\cite{We95a}.

\subsection{The case $\Phi(M,t_0,t_\mathbf{f})=M$}
\label{subsec:first_case}

We first choose a set $Y=\{y_1,\ldots,y_n\}\subset M$ of \emph{centers}.  The corresponding basis functions $\varphi_j:M\to\R$ are given by
\[
\varphi_j(x)=\psi(\eps\|x-y_j\|_2),
\]
$j=1,\ldots,n$, where $\eps>0$ is the \emph{shape parameter} which scales the size of the support.

The eigenfunctions $f$ of $\triangle^D$ will be approximated by functions from $\cA=\spn(\varphi_1,\ldots,\varphi_n)\subset L^2(M)$, i.e.\ by linear combinations of the form
\[
\tilde f(x) = \sum_{j=1}^n \alpha_j\varphi_j(x).
\]
To this end, we choose $\nin$ inner collocation nodes $X_\inn =\{x_1,\ldots,x_{\nin}\}\subset \intr M$ in the interior of $M$ as well as $\nbd=n-\nin$ nodes $X_\bd =\{x_{\nin+1},\ldots,x_n\}$ on the boundary of $M$.  If $M$ has no boundary, then $\nin=n$ and $X_\bd =\emptyset$.
 Now, given a (sufficiently smooth) function $f:M\to\R$ which satisfies the boundary condition
\[
\cL_\bd f:=\nabla f(x)\cdot\left[(I+C_{x,t_0,t_\mathbf{f}}^{-1})\mathbf{n}(x)\right]=0
\] on $\partial M$, the coefficients $\alpha=(\alpha_1,\ldots,\alpha_n)^T$ of its interpolating approximation $\tilde f\in\cA$ are given by the solution of the linear system
\[
A\alpha = E_0 f_\inn
\]
where
\[
A = \begin{bmatrix}A_\inn \\ L_\bd \end{bmatrix}, \quad
E_0 = \begin{bmatrix}I\\ 0\end{bmatrix}, \quad
f_\inn =\begin{bmatrix}f(x_1)\\\vdots\\f(x_{\nin})\end{bmatrix},
\]
and
\begin{align*}
A_\inn  & =  \left(\varphi_j(x_i)\right)_{ij}, &  i & =1,\ldots,\nin, & j & =1,\ldots,n\\
L_\bd  & = \left(\cL_\bd \varphi_j(x_i)\right)_{ij}, & i & =\ell+1,\ldots, n, & j & = 1,\ldots,n.
\end{align*}

For the  case $\partial M=\emptyset$, the fact that the interpolation matrix $A_\inn$ is (strictly) positive definite implies the invertibility of $A$.  In the case $\partial M\neq\emptyset$, however, invertibility does not seem to be guaranteed for an arbitrary choice of the nodes $X_\inn, X_\bd$ \cite{PlDr03a}. In our numerical experiments, however, $A$ was always invertible.

Given some arbitrary linear operator $\cL:L^{2}(M)\to L^{2}(M)$, the image $\cL\tilde f$ of a function $\tilde f\in\cA$ is given by
\[
\cL\tilde f(x) = \sum_{j=1}^n \alpha_j \cL\varphi_j(x).
\]
Correspondingly, the values $\tilde g_i = \cL\tilde f(x_i)$ of this image at the inner collocation nodes $x_i$, $i=1,\ldots,\nin$, can be computed by
\[
\tilde g = L_\inn  \alpha = L_\inn  A^{-1}E_0 \tilde f_\inn ,
\]
where
\[
L_\inn  =  \left(\cL\varphi_j(x_i)\right)_{ij}, \quad i=1,\ldots,\nin, \quad j=1,\ldots,n.
\]
That is, the matrix $L:=L_\inn  A^{-1}E_0$ can be seen as a discrete version of the given linear operator $\cL$.

In our case, we are concerned with the linear operators $\cP$, $\triangle$ and $\cP^*$ (cf.\ \ref{hattriangle}) and so the approximating operators are $P=P_\inn  A^{-1}E_0$ with
\[
P_\inn  =  \left(\cP\varphi_j(x_i)\right)_{ij}, \quad i=1,\ldots,\nin, \quad j=1,\ldots,n,
\]
$P^*=P^*_\inn  A^{-1}E_0$ with
\[
P^*_\inn  =  \left(\cP^*\varphi_j(x_i)\right)_{ij}, \quad i=1,\ldots,\nin, \quad j=1,\ldots,n,
\]
and $D=D_\inn  A^{-1}E_0$ with
\[
D_\inn  =  \left(\triangle\varphi_j(x_i)\right)_{ij}, \quad i=1,\ldots,\nin, \quad j=1,\ldots,n.
\]
The discrete eigenproblem associated to (\ref{hattriangle})-(\ref{strongbc0}) is
\[
\frac12(D+P^*DP)\tilde f = \lambda \tilde f.
\]
We note that $\mathcal{P}$ is unitary, so $\mathcal{P}^*=\mathcal{P}^{-1}$.  In general, this property is not inherited by the matrix $P$ which is constructed via collocating the action of $T^{-1}$.  Likewise, the matrix $P^*$ is constructed by collocating the action of $T$, i.e.\ one may view $P^*$ as taking the dual of $\mathcal{P}$ followed by collocation.  In general, this approach will not guarantee that $P^*$ is the dual of $P$.  As a result, the discrete dynamic Laplacian $\frac12(D+P^*DP)$ is not necessarily self-adjoint and therefore its spectrum not necessarily real (in contrast to $\triangle^D$)\footnote{In addition, the matrix $D$ is not necessarily self-adjoint.}.

However, if we interchange these two operations, i.e.\ if we first collocate and then take the dual, we obtain $P^\top$ which is the dual of $P$.  In our numerical experiments, using $P^\top$ instead of $P^*$ turned out to yield a much improved spectrum, i.e.\ in Section~\ref{sec:experiments} we solve the eigenproblem
\begin{equation}\label{eq:discrete_eigenproblem}
\frac12(D+P^\top DP)\tilde f = \lambda \tilde f.
\end{equation}

\textbf{Remark.} Note that, as mentioned above, this discretization  applies to an arbitrary linear operator $\cL$ (cf.\  Platte and Driscoll \cite{PlDr03a}).  In particular, one may apply this to the ``classical'' approach for computing coherent pairs \cite{F13}, namely discretizing the transfer operator of a small random perturbation of the underlying system.  In this case, in addition to the matrix $P$, one would compute a discretization $D^\eps$ of some local averaging operator $\mathcal{D}^\eps$ and compute the largest singular values and vectors of the matrix $D^\eps P D^\eps$.

\subsection{Implementational details}
\label{sect:implementation}

In implementing the approach described above, various choices on the numerical parameters have to be made which strongly influence the resulting approximation quality.  Here, we collect a few comments on a concrete implementation. Corresponding Matlab codes can be obtained from the authors.

We first need to decide on a set $Y$ of centers for the basis functions.  While in general, scattered points are fine for approximations with radial basis functions (in fact, this is often one of the main motivations for using radial basis functions), in our experiments regular grids led to a much improved accuracy.  In our 2D examples, a set $Y$ of 1000-2000 centers led to sufficiently accurate results.

The next choice concerns the shape parameter $\eps$ which governs the radius of the support of the basis functions.  Again, the computational results are typically highly sensitive with respect to this parameter.  In general, better results are obtained when choosing a smaller value for $\eps$, leading to larger supports.  On the other hand, the condition number of the interpolation matrices (i.e.\ $A_\inn $, etc.) increases with decreasing $\eps$ and for values too small, instabilities will occur.  In our experiments, we chose $\eps$ such that the support of a basis function overlapped with roughly 100-200 other ones.

For the collocation points $X_\inn , X_\bd$, the same remarks apply as for the centers.  Here, we use choose $X_\inn\cup X_\bd=Y$ since this seems to preserve the spectral properties of the various matrices best.  Only in the second example (cf.\ Section~\ref{sec:experiments}), we slightly shifted the centers away from the boundary in order to keep the code simpler.  Note that the choice $X_\inn\cup X_\bd=Y$ requires a special treatment of those entries of the discrete Laplacian where the collocation point coincides with the center (because of the zero in the denominator of the Laplacian in polar coordinates, which we employ due to the radial structure of our basis functions).

To compute the inverse of the Cauchy-Green tensor, we integrate the variational equation backward in time, i.e.\ explicitly compute $D_x\Phi(\Phi(x,t_0,t_\mathbf{f}),t_\mathbf{f},t_0)$ for $x\in X_\bd $ backwards along the nonlinear trajectory $\{\Phi(x,t_0,t)\}_{t_0\le t\le t_\mathbf{f}}$ computed prior in forward time.  For longer flow times, this may yield a (nearly) singular matrix.

We also have to compute the various matrices: For $A_\bd $ and $L_\inn $, we first apply the differential operators to the basis functions $x\mapsto \varphi_j(x)$ and then evaluate at the collocation points. Similarly, when computing $P_\inn $, we compute
\[
(P_\inn )_{ij} = \varphi_j(\Phi(x_i,t_\mathbf{f},t_0)).
\]

Finally, for the resulting discrete eigenproblem (\ref{eq:discrete_eigenproblem}), we compute a few eigenvalues of smallest magnitude by using an implicitly restarted Arnoldi method \cite{LeSoYa98a} as implemented in the \texttt{eigs} function in Matlab (i.e.\ with parameter \texttt{'SM'}).

\subsection{The general case $\Phi(M,t_0,t_\mathbf{f})\neq M$}

Whenever $\hat M := \Phi(M,t_0,t_\mathbf{f})\not\subset M$, it will be necessary to work with two distinct approximation spaces that are supported on $M$ and $\hat M$, respectively.
In the context of radial basis function approximations, it appears natural to define the set of centers $\hat Y$ in $\hat M$ by the images of the centers in $Y$ under the flow $\Phi$ (in this section for brevity we will write $\Phi$ for $\Phi(\cdot,t_0,t_\mathbf{f})$) which here we assume to be a bijection (as is the case if the underlying vector field $F$ is locally Lipschitz w.r.t.\ $x$).  In fact, the discrete dynamic Laplacian attains a particular simple form in this case as we will see now.  Choosing $\hat Y = \{\Phi(y) \mid y\in Y\}$,  we have the associated basis functions
\[
\hat\varphi_j(x)=\psi(\eps\|x-\hat y_j\|_2),
\]
$j=1,\ldots,n$, so that we can approximate functions $\hat f:\hat M\to\R$ by elements from the space $\hat\cA=\spn(\hat\varphi_1,\ldots,\hat\varphi_n)$.  Similarly, we define the set of collocation points as $\hat X_\inn =\Phi(X_\inn )=\{\hat x_1,\ldots,\hat x_n\}=\{\Phi(x_1),\ldots,\Phi(x_n)\}$ and $\hat X_\bd  = \Phi(X_\bd )$.

For the case $\partial M = \emptyset$ one can proceed as follows:  Given some function $f(x)=\sum_{j=1}^n \alpha_j\varphi_j(x)\in\cA$, the approximation of $\cP f$ within $\hat\cA$ will be given by a linear combination of the form $\sum_{j=1}^n \hat\alpha_j\hat\varphi_j(x)$.  In fact, requiring this equality to hold in the collocation points $\hat x_i$, we obtain the system
\[
\cP f(\hat x_i) = \sum_{j=1}^n \alpha_j \varphi_j(\Phi^{-1}(\hat x_i)) = \sum_{j=1}^n \hat\alpha_j\hat\varphi_j(\hat x_i),
\]
$i=1,\ldots,n$.  Since $\Phi^{-1}(\hat x_i)=\Phi^{-1}(\Phi(x_i))=x_i$, we obtain the matrix representation
\[
P_\alpha = \hat A^{-1}A,
\]
for the discretization of $\cP$ as a map (on the coefficients $\alpha$) from $\cA$ to $\hat\cA$, where $A_{ij}=(\varphi_j(x_i))$ and $\hat A_{ij}=(\hat\varphi_j(\hat x_i))$, $i,j=1,\ldots,n$.  Similarly, for $\cP^*$, we get the matrix
\[
P^*_\alpha = A^{-1}\hat A.
\]
Again, one might want to use $P^\top_\alpha$ instead of $P^*_\alpha$ in order to enforce self adjointness of the discrete dynamic Laplacian, cf.\ the comment at the end of Section~\ref{subsec:first_case}.

For the application of the Laplacian $\triangle$, we again have to distinguish whether it is being applied to a function from $\cA$ or $\hat\cA$.  Mimicking the above derivation for $P$, we obtain the matrix representation
\[
D_\alpha = A^{-1}D'
\]
for the discretization of $\triangle$ as a map (again, on the coefficients $\alpha$) from $\cA$ to $\cA$, where $D' = (\triangle\varphi_j(x)|_{x=x_i})_{ij}$.  Likewise, for $\hat \triangle$ (we call the second $\triangle$ in (\ref{hattriangle}) $\hat{\triangle}$ as it acts on $\Phi(M)$)
\[
\hat D_\alpha=\hat A^{-1}\hat D'
\]
with $\hat D' = (\triangle\hat\varphi_j(x)|_{x=\hat x_i})_{ij}$ and the discrete dynamic Laplacian is
\begin{equation}
\label{coeffdynlapeqn}
\frac12(D_\alpha+P^*_\alpha\hat D_\alpha P_\alpha).
\end{equation}

Note that the matrix (\ref{coeffdynlapeqn}) is a mapping on the coefficient vectors $\alpha$ of some RBF-approximation $\tilde f\in\cA$ -- in contrast to (\ref{eq:discrete_eigenproblem}), where the matrices are mappings on vectors $f_\inn \in\R^{\nin}$ of function values.  Certainly, we can rephrase the eigenproblem here in terms of $f_\inn $, too.  To this end, we apply $A$ and $A^{-1}$ (resp.\ $\hat A$ and $\hat A^{-1}$) appropriately, yielding the matrices
\begin{align*}
P & = \hat A\hat A^{-1}A A^{-1} = I \\
P^* & = A A^{-1}\hat A \hat A^{-1}  = I \\
D & = AA^{-1}D'A^{-1} = D'A^{-1} \\
\hat D & = \hat A\hat A^{-1}\hat D'\hat A^{-1} = \hat D'\hat A^{-1},
\end{align*}
i.e.\ we obtain the discrete dynamic Laplacian
\begin{equation}
\label{eq:dis_dyn_Lap_f}
\frac12(D+\hat D).
\end{equation}
Let us explore a little the similarity between the expressions (\ref{hattriangle}) and (\ref{eq:dis_dyn_Lap_f}).
In the second term of (\ref{eq:dis_dyn_Lap_f}) we have $(\hat D')_{ij}=\triangle\hat\varphi_j(x)|_{x=\hat x_i}=\triangle\psi_j(\|x-\Phi(y_j)\|_2)|_{x=\Phi(x_i)}$.
That is, in computing $\hat D'$, we are using the distance matrix $(\|\Phi(x_i)-\Phi(y_j)\|_2)_{ij}$ of pairwise distances between the \emph{images} of the centers and the collocation points.
This is consistent with the fact that the operator $\mathcal{P}^*\triangle\mathcal{P}$ is the Laplace operator on $M$ endowed with the $\Phi$-pullback of the Euclidean metric on $\Phi(M)$; see Section 3.2 \cite{F15} or p27 \cite{chaveleigenvalues}.

\subsection{Extraction of coherent sets from trajectories}
\label{sect:extract}

The algorithm we use in the following two-dimensional volume-preserving case studies  is described below.  We consider the case $\Phi(M,t_0,t_\mathbf{f})=M$ and correspondingly use the approach described in Section \ref{subsec:first_case}.
Algorithm 1 is presented for the situation where the boundary size is compared at the initial and final times;  the obvious extensions can be made if additional comparison times in the interval $[t_0,t_\mathbf{f}]$ are desired.  We outline the main steps here, more detailed comments on the individual steps can be found in Section~\ref{sect:implementation}.
\vspace{.2cm}

\noindent\textbf{Algorithm 1:}
\begin{enumerate}
\item Choose a set of centres $Y=\{y_j\}_{j=1}^n$ in $M$ and a shape parameter $\epsilon > 0$.
\item Select a set of internal collocation points $X_\inn =\{x_i\}_{i=1}^{\nin}$ and a set of boundary collocation points $X_\bd =\{x_i\}_{i=\nin+1}^n$ in $\partial M$.  If $M$ has no boundary, then $\nin=n$ and $X_\bd =\emptyset$.
\item If $X_\bd \neq\emptyset$, calculate the unit normal $\mathbf{n}(x_i)$ and the Cauchy-Green tensor $C_{x_i,t_0,t_\mathbf{f}}$ for $i=\OJ{\nin}+1,\ldots,n$.
\item  Form the matrix ${P}$ for the transfer operator and the matrix $D$ for the Laplacian as described in Section \ref{subsec:first_case} as well as $D'=(D+P^\top DP)/2$.
\item Solve the eigenproblem (\ref{eq:discrete_eigenproblem}), \OJ{i.e.\ compute several eigenvalues $\lambda_1>\lambda_2>\cdots$ of smallest magnitude (e.g.\ using Matlab's \texttt{eigs} command with the option \texttt{'SM'})} and \OJ{associated} eigenvectors $f_1, f_2,\ldots$ of $D'$.
\item Iteratively scan over values of $f_2$ from $\min f_{2}$ to $\max f_{2}$. For each value, extract a level curve $\Gamma$ in $M$ using Matlab's \verb"contour" function (this function returns a collection of points representing corners of a polygonal curve).   To compute $\Phi(\Gamma,t_0,t_\mathbf{f})$,
\begin{enumerate}
\item Either: map the points representing $\Gamma$ directly with $\Phi(\cdot,t_0,t_\mathbf{f})$,
\item Or: compute $Pf_2$ and extract $\Phi(\Gamma,t_0,t_\mathbf{f})$ using Matlab's \verb"contour" function with the same level set value as for $\Gamma$.
\end{enumerate}
\item Optimise $\mathbf{h}^D(\Gamma)$ by running over
all curves $\Gamma$ formed from level sets of $f_2$ in Step 6.
The length of $\Gamma$ and $\Phi(\Gamma,t_0,t_\mathbf{f})$ are computed as the lengths of the polygonal curves comprising them. Report the $\Gamma$ and $\Phi(\Gamma,t_0,t_\mathbf{f})$ that yield the lowest value of $\mathbf{h}^D(\Gamma)$.
\end{enumerate}

\section{Numerical experiments}
\label{sec:experiments}

\subsection{Standard map on the torus}
\label{sec:standard}

\begin{figure}[htbp]
  \centering
	\begin{tikzpicture}
    \begin{axis}[enlargelimits=false,width=0.33\textwidth,height=0.33\textwidth,
		xtick={0,6.28},ytick={0,6.28},xticklabels={0,$2\pi$},yticklabels={\empty,$2\pi$},
		xlabel={$x$},xlabel style={at={(0.5,0.1)}},ylabel={$y$},ylabel style={at={(0.2,0.5)}}]
        \addplot graphics [xmin=0,xmax=6.28,ymin=0,ymax=6.28]{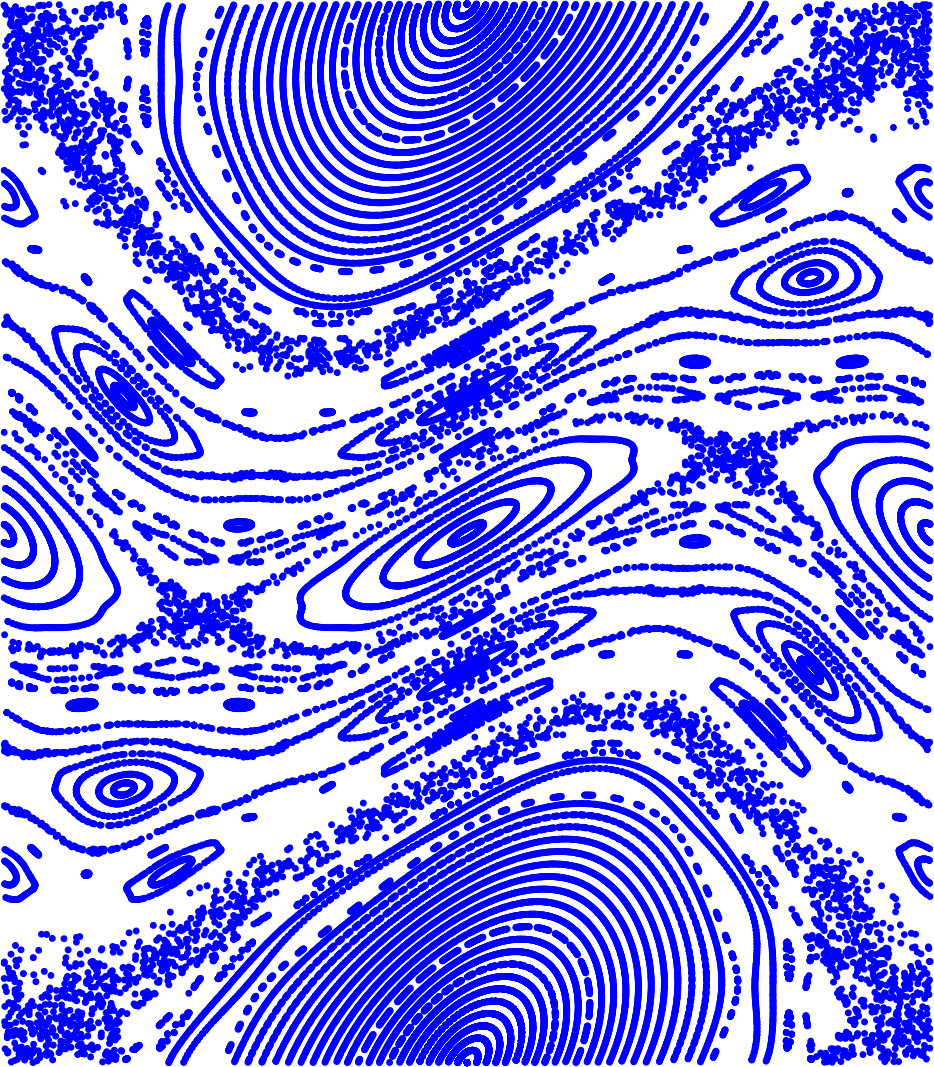};
    \end{axis}
	\end{tikzpicture}
  \caption{Orbits of the map (\ref{stdmap}) showing regular and chaotic regions.}\label{fig:infaction}
\end{figure}

\begin{figure}[htbp]
  \centering
	\begin{tikzpicture}
    \begin{axis}[enlargelimits=false,width=0.25\textwidth,height=0.25\textwidth,
		xtick={0,6.28},ytick={0,6.28},xticklabels={0,$2\pi$},yticklabels={\empty,$2\pi$},
		xlabel={$x$},xlabel style={at={(0.5,0.1)}},ylabel={$y$},ylabel style={at={(0.3,0.5)}}]
        \addplot graphics [xmin=0,xmax=6.28,ymin=0,ymax=6.28]{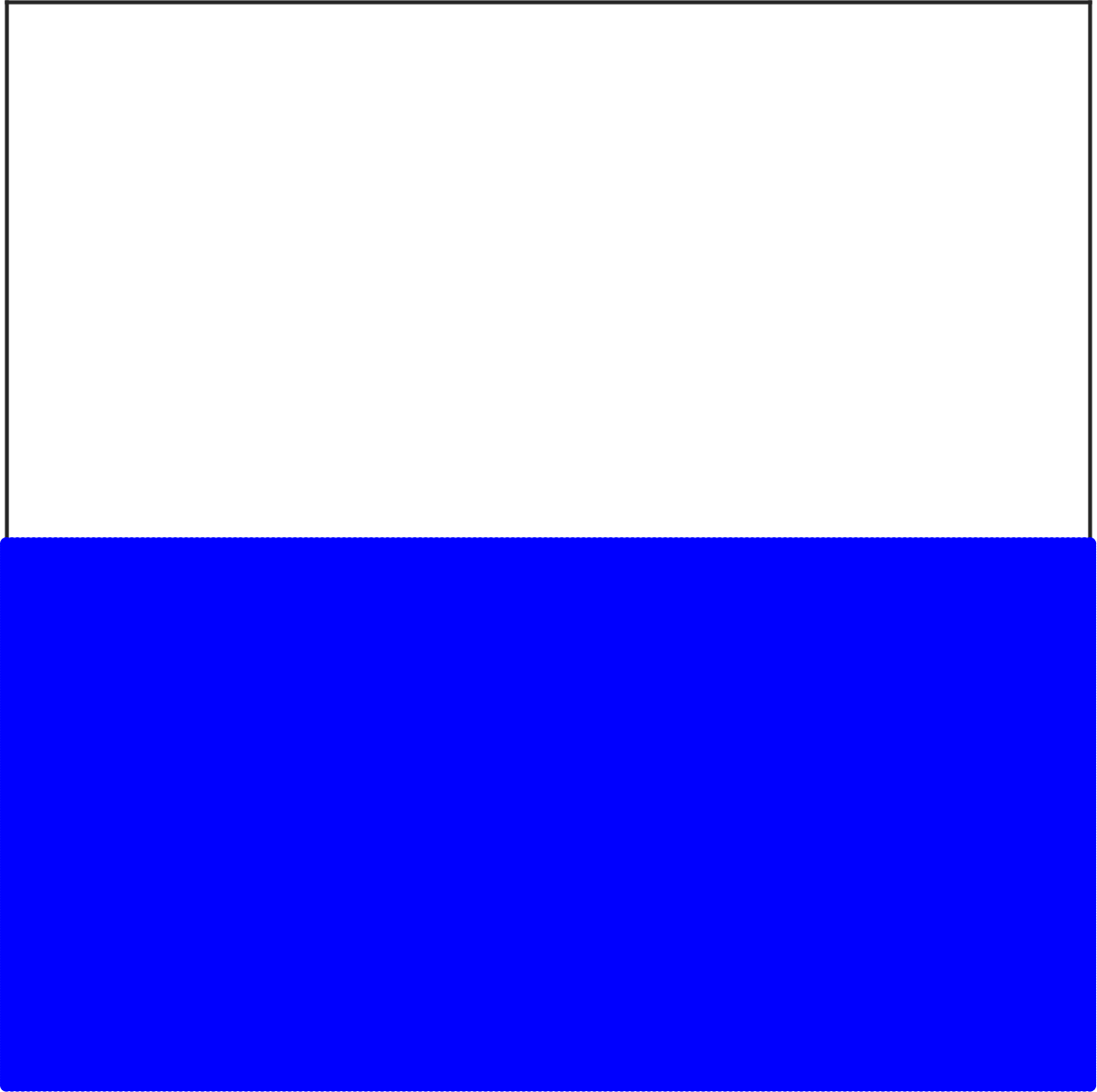};
    \end{axis}
	\end{tikzpicture}
	\begin{tikzpicture}
    \begin{axis}[enlargelimits=false,width=0.25\textwidth,height=0.25\textwidth,
		xtick={0,6.28},ytick={0,6.28},xticklabels={0,$2\pi$},yticklabels={\empty,$2\pi$},
		xlabel={$x$},xlabel style={at={(0.5,0.1)}},ylabel={$y$},ylabel style={at={(0.3,0.5)}}]
        \addplot graphics [xmin=0,xmax=6.28,ymin=0,ymax=6.28]{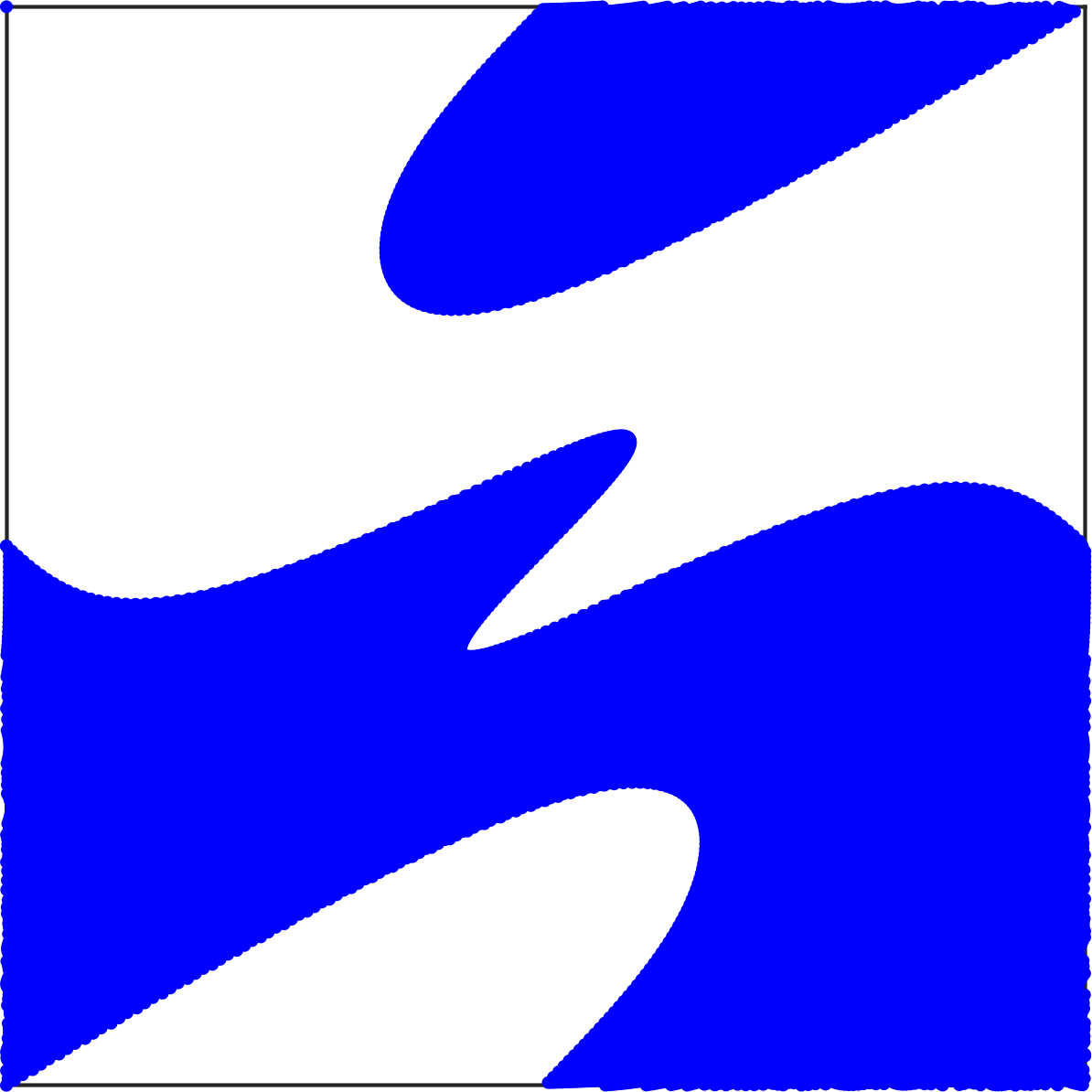};
    \end{axis}
	\end{tikzpicture}
  \caption{Left: An initial colouring of phase space. Right: The image of the colouring under $T^2$.}\label{fig:action}
\end{figure}

We consider the standard map
\begin{equation}
\label{stdmap}
T(x,y)=(x+y+a\sin x,y+a\sin x)\pmod{2\pi},
\end{equation}
with parameter $a=0.971635$, a value above which a prominent KAM curve is destroyed, and the phase space shows both regular and chaotic motion; see Fig.~\ref{fig:infaction}.
We are interested in finite-time dynamics rather than time-asymptotic dynamics and we choose to analyse the dynamics of two iterates of $T$.
The geometric action of $T^2$ in phase space is indicated in Figure \ref{fig:action}.

We applied Algorithm 1 as follows: 1. We choose the centers $Y = h\Z^2\cap [0,2\pi]^2$ on a regular grid with $h=0.33$, leading to $400$ points in $Y$,  use the Wendland function $\psi_{6,4}$ for the basis functions and $\eps=0.4$ for the value of the shape parameter, such that the support of each basis function overlaps with roughly 200 other basis functions.  2. We choose $X_\inn = Y$ for the inner collocation points and since $\partial M=\emptyset$ we have $X_\bd =\emptyset$ and we do not need to compute any normals or Cauchy-Green tensors.   3. We compute the matrices $P$ and $D$ as described in Sections~\ref{subsec:first_case} and \ref{sect:implementation} and 4. solve the eigenproblem (\ref{eq:discrete_eigenproblem}).

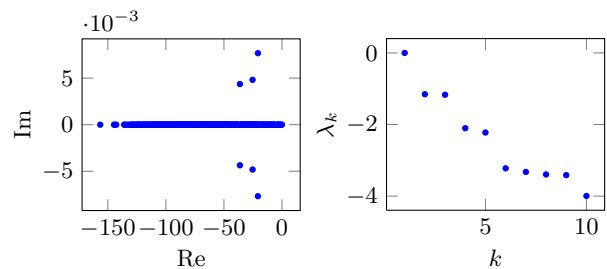
\begin{figure}[htbp]
  \centering
   \begin{tikzpicture}
    \begin{axis}[
    	width=0.25\textwidth,
    	ylabel={Im},ylabel style={at={(0.15,0.5)}},
    	xlabel={Re},xlabel style={at={(0.5,0.05)}},
		mark size = 1
    ]
    \addplot[blue,only marks] table {stdmap3_specDL.tex};
    \end{axis}
  \end{tikzpicture}
   \begin{tikzpicture}
    \begin{axis}[
    	width=0.25\textwidth,
    	ylabel={$\lambda_k$},ylabel style={at={(0.15,0.5)}},
    	xlabel={$k$},xlabel style={at={(0.5,0.05)}},
		mark size = 1
    ]
    \addplot[blue,only marks] table {stdmap3_specDL_close.tex};
    \end{axis}
  \end{tikzpicture}

  \caption{Standard map: Spectrum of $\frac12 (D+P^\top DP)$ (left: entire spectrum, right: eigenvalues closest to~$0$).}\label{fig:std_spectrum}
\end{figure}

The resulting four leading eigenvalues were (to three significant figures):   $\lambda_1=6\cdot 10^{-5}$, $\lambda_2=-1.15$,   $\lambda_3=-1.17$ and $\lambda_4=-2.10$ which appear to be correct when compared to the results of a highly accurate computation using a spectral method.  Fig.~\ref{fig:std_spectrum} shows the spectrum of the discrete dynamic Laplacian $\frac12 (D+P^\top DP)$.
Note that this spectrum is not real (as it should be).
However, since the part of the eigenspectrum we are interested in is real, namely in a neighborhood of the origin (see Fig.~\ref{fig:std_spectrum} (right)), we do not need to employ further techniques to reduce the imaginary components of the spectrum of $\frac12 (D+P^\top DP)$ (such as using a least squares approach instead of interpolation\cite{PlDr06a}).
One could also enforce a real spectrum by using $(D+D^\top)/2$ instead of $D$.
The question of how these two techniques change the spectrum is currently under investigation.

Note that in order to obtain this accuracy, we only evaluated the map $20\cdot 20=400$ times. \OJ{Experimentally, in order to obtain the same accuracy with Ulam's method one needs to use a grid of $64\times 64=4096$ boxes and $7\times 7=49$ sampling points per box at least, requiring the mapping of $4096\cdot 49 \approx 2\cdot 10^5$  points.  On the other hand for our current approach, setting up the distance matrices as described in Section \ref{subsec:first_case} requires additional computational effort and the resulting eigenproblem is not sparse, in contrast to the one resulting from Ulam's method.}

\begin{figure}[hbt]
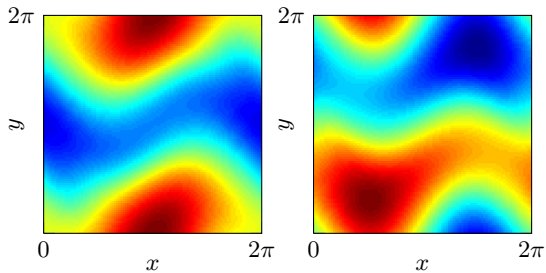

  \hspace*{-8mm}
  \evplot{std_map_ev2}{}{$x$}{$y$}\hspace*{-4mm}
  \evplot{std_map_ev3}{}{$x$}{$y$}
  \caption{Two iterations of the standard map: Eigenvectors $f_2$ (left) and $f_3$ (right) of (\ref{eq:discrete_eigenproblem}), using $n=400$ centers on a regular grid with $\eps=0.4$.}\label{fig:fineevecs}
\end{figure}

\begin{figure}[hbt]
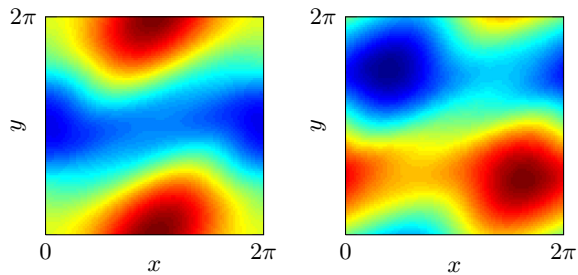

  \centering
  \evplot{std_map_ev2_img}{}{$x$}{$y$}
  \evplot{std_map_ev3_img}{}{$x$}{$y$}
  \caption{Images $Pf_2$ (left) and $Pf_3$ (right) of the eigenvectors from Figure \ref{fig:fineevecs}.}\label{fig:finepushevecs}
\end{figure}

The eigenvalues $\lambda_2$ and $\lambda_3$ are of similar value, indicating that there are two independent ways to disconnect $M$ such that the disconnections each have similarly small dynamic boundary size;  see Figs.\ \ref{fig:fineevecs} and \ref{fig:finepushevecs}.
The first two nontrivial eigenvectors are shown in Figure  \ref{fig:fineevecs}. {In producing this figure (as well as all subsequent figures of eigenvectors of the dynamic Laplacian), we evaluated the RBF representation of the eigenvectors on a grid of $100\times 100$ points.}
Their images $Pf_2$ and $Pf_3$ under $P$ are shown in Figure \ref{fig:finepushevecs};  these vectors correspond to the geometry of the phase space \emph{after} two iterations of $T$.

\begin{figure}[hbt]
  \centering

   \begin{tikzpicture}
    \begin{semilogyaxis}[
    	width=0.4\textwidth,
    	ylabel={$h^D_{\{0,2\}}(\Gamma_\gamma)$},
    	xlabel={$\gamma$},
    	grid
    ]
    \addplot[blue] table {std_map_cheeger.tex};
    \end{semilogyaxis}

  \end{tikzpicture}

  \caption{Standard map: Graph of the dynamic Cheeger constant (\ref{cheegereqn}) in dependence of the level set value $\gamma$.}\label{fig:thresh}
\end{figure}
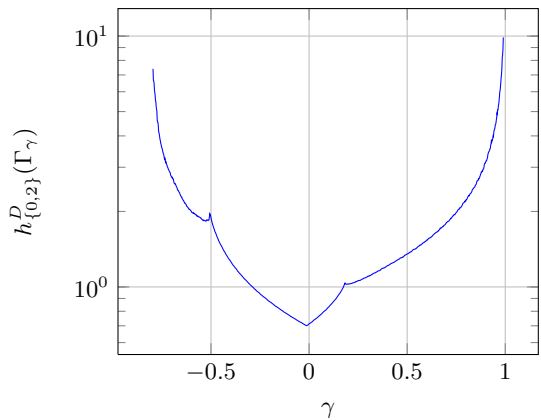

Step 6 of Algorithm 1 involves scanning over the level sets of $f_2$ (shown in the left in Figures \ref{fig:fineevecs} and \ref{fig:finepushevecs}) and calculating the length of the level set and the length of its image.
We take the approach of version (b) in Step 6, and Figure \ref{fig:thresh} shows the value of (\ref{cheegereqn}) as a function of $\gamma\in\mathbb{R}$, where $\Gamma_\gamma=\{x\in \mathbb{T}^2: f_2(x)=\gamma\}$.

\begin{figure}[hbt]
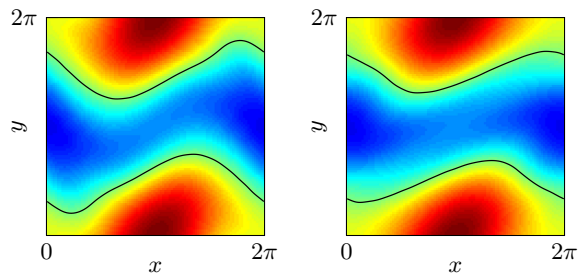

  \centering
  \evplot{stdmap3_ev2_curve}{}{$x$}{$y$}
  \evplot{stdmap3_ev2_curve_img}{}{$x$}{$y$}

  \caption{Minimising curves before (left) and after (right) the application of $T^2$.}\label{fig:curves}
\end{figure}

\begin{figure}[hbt]
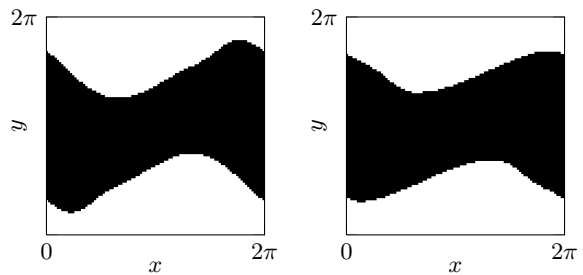

  \centering
  \evplot{stdmap3_ev2_cs}{}{$x$}{$y$}
  \evplot{stdmap3_ev2_cs_img}{}{$x$}{$y$}
  \caption{Finite-time coherent sets before (left) and after (right) the application of $T^2$.}\label{fig:sets}
\end{figure}

The minimising $\gamma=-0.0115$ was chosen by scanning through 1000 equally separated values of $\gamma$ from $\min f_2$ to $\max f_2$, where $f_2$ is scaled so that $|f_2|_\infty=1$.
The curves along the corresponding level sets in $f_2$ and $Pf_2$ are shown in Figure \ref{fig:curves}, having lengths $\ell_1(\Gamma_\gamma)=14.90$ and $\ell_1(\Phi(\Gamma_\gamma,t_0,t_\mathbf{f}))=13.54$, while $\ell_2(M_1)=20.12$ for the corresponding finite-time coherent sets $M_1$ shown in Figure \ref{fig:sets}, so that $h^D_{0,2}(\Gamma_\gamma)=0.70$.
Note that the curves in Figure \ref{fig:curves} (left) and (right) are similar, but not identical.
There is no prescription for the curves to be invariant under $T^2$;  the fact that they are approximately invariant is likely tied to the asymptotic dynamics of $T$.

The curves in Figure \ref{fig:curves} are in the vicinity of the ``inner'' boundaries (the boundaries closer to $y=0.5$) of the chaotic regions in Figure \ref{fig:infaction}.
The {red} features in Figure \ref{fig:curves} also correspond reasonably closely to the largest regular regions in the upper and lower parts of Figure \ref{fig:infaction}.
Figure \ref{fig:infaction} displays the time-asymptotic dynamics of $T$, while we are analysing two iterations of $T$.
If we increased the number of iterations in our analysis, it is likely that the level sets of $f_2$ will more closely reflect the time-asymptotic dynamics.

The third eigenfunction $f_3$ shows further finite-time structures that undergo low levels of deformation (the boundaries of the approximately red and blue sets remain small at the initial (Figure \ref{fig:fineevecs} (right)) and final (Figure \ref{fig:finepushevecs} (right)) times, respectively.
The total boundary lengths of the red/blue interfaces in Figures \ref{fig:fineevecs} and \ref{fig:finepushevecs} are larger for $f_3$ than for $f_2$ (albeit only slightly so).
These finite-time structures are highly mobile from time $t=0$ to time $t=2$ and are not at all evident in Figure \ref{fig:infaction}.
	
\paragraph*{Order of convergence.} One of the benefits of using radial basis functions for function approximation is that, depending on the ``basic'' function $\psi$ one chooses, high convergence orders can be obtained.  Here, we have chosen the compactly supported Wendland functions $\psi=\psi_{s,k}$.

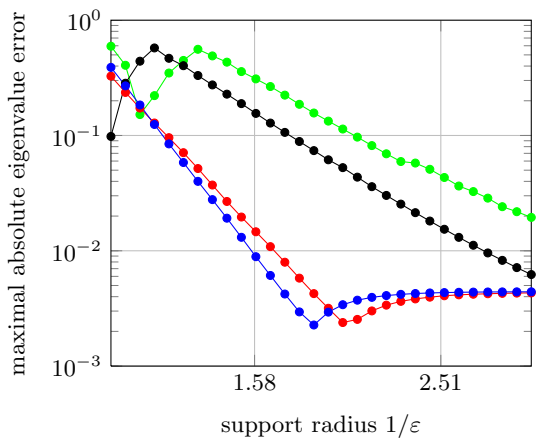
\begin{figure}[htbp]
\centering
\begin{tikzpicture}
\begin{loglogaxis}[
	enlargelimits=false,
	width=0.4\textwidth,
	xlabel={support radius $1/\varepsilon$},
	ylabel={maximal absolute eigenvalue error},
	ymin=1e-3,
	ymax=1e0,
	xticklabel=\pgfmathparse{exp(\tick)}\pgfmathprintnumber{\pgfmathresult},
	mark size=1.5pt,
	grid
]
\addplot[green,mark=*] table {stdmap_ev_err_vs_eps_N20_phi31.tex};
\addplot[black,mark=*] table {stdmap_ev_err_vs_eps_N20_phi42.tex};
\addplot[red,mark=*] table {stdmap_ev_err_vs_eps_N20_phi53.tex};
\addplot[blue,mark=*] table {stdmap_ev_err_vs_eps_N20_phi64.tex};
\end{loglogaxis}
\end{tikzpicture}

	\caption{Standard map: Dependence of the maximal absolute error in the eigenvalues $\lambda_1,\ldots,\lambda_4$ on  the support radius of the basis functions for $n=400$ centers on an equidistant grid. Green: $\psi_{3,1}$, black: $\psi_{4,2}$, red: $\psi_{5,3}$, blue: $\psi_{6,4}$).}\label{fig:std_eval_error_vs_eps_1}
\end{figure}

We first investigate how the absolute errors in the largest few eigenvalues depend on the shape parameter.  Figure~\ref{fig:std_eval_error_vs_eps_1} shows the results of an  experiment for $n=400$ centers.  Here, the results of the RBF approach are compared to those of a computation using spectral collocation which appears to be highly accurate.
In this figure, the maximal error in the four eigenvalues with smallest magnitude in dependence of the radius $1/\eps$ of the support of the basis functions $\phi_j$ is shown for various choices of the basic function $\psi_{s,k}$.  Clearly, these errors sensitively depend on $\eps$ and smoother basis functions yield lower errors.  With $\psi=\psi_{6,4}$ we obtain an error of approximately $10^{-3}$ for $1/\eps\approx 1.6$. Recall that we only needed to evaluate the map $n=400$ times here.  On the other hand, for $\psi=\psi_{6,4}$ and $\psi=\psi_{5,3}$, a further decrease in $\eps$ does not yield a smaller error.

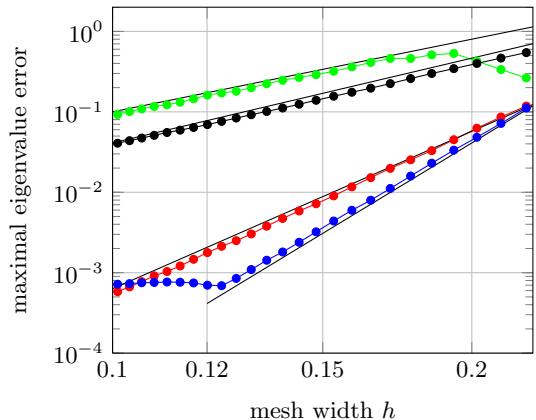
\begin{figure}[htbp]
\centering
\begin{tikzpicture}
\begin{loglogaxis}[
	enlargelimits=false,
	width=0.4\textwidth,
	xlabel={mesh width $h$} ,
	ylabel={maximal eigenvalue error},
	ymin=1e-4,
	ymax=2,
	xticklabel=\pgfmathparse{exp(\tick)}\pgfmathprintnumber{\pgfmathresult},
	xtick={0.1,0.12,0.15,0.2},
	mark size=1.5pt,
	grid
]
\addplot[green,mark=*] table {stdmap_ev_err_vs_h_eps08_phi31.tex};
\addplot[black,mark=*] table {stdmap_ev_err_vs_h_eps08_phi42.tex};
\addplot[red,mark=*] table {stdmap_ev_err_vs_h_eps08_phi53.tex};
\addplot[blue,mark=*] table {stdmap_ev_err_vs_h_eps08_phi64.tex};
\addplot[black,domain=0.1:0.225,samples=2] {1e2*x^3};
\addplot[black,domain=0.1:0.225,samples=2] {1.3e2*x^3.5};
\addplot[black,domain=0.1:0.225,samples=2] {2e3*x^6.5};
\addplot[black,domain=0.12:0.225,samples=2] {8e4*x^9};
\end{loglogaxis}
\end{tikzpicture}

	\caption{Standard map:  Dependence of the maximal absolute error in the first 4 eigenvalues on the mesh widths for $\eps=0.8$. Green: $\psi_{3,1}$, black: $\psi_{4,2}$, red: $\psi_{5,3}$, blue: $\psi_{6,4}$).
}\label{fig:std_eval_error_vs_h_1}
\end{figure}

We next analyse how the eigenvalue error depends on the number of centers, i.e.\ the mesh width.  To this end, we need to fix the value of the shape parameter $\eps$ since it is known that an interpolating approximation with RBFs in which one scales the support of the basis functions proportionally to the mesh width does not converge\cite{We05a,Fa07a}.  We perform an experiment with $\eps=0.8$ and show the maximal error in the four eigenvalues with smallest magnitude in Figure~\ref{fig:std_eval_error_vs_h_1}.  We observe convergence orders of $\cO(\eps^{-3}), \cO(\eps^{-3.5}),\cO(\eps^{-6.5})$ and $\cO(\eps^{-9})$ as indicated by the black lines (from top to bottom) by using $\psi_{3,1}$, $\psi_{4,2}$, $\psi_{5,3}$ and $\psi_{6,4}$, respectively.  Clearly, using basis functions of higher smoothness pays off, but surprisingly, $\psi_{3,1}$ and $\psi_{4,2}$ seem to converge with approximately the same rate (albeit $\psi_{4,2}$ delivers a smaller error).  Finally, in Figure~\ref{fig:std_eval_error_vs_h_1}, the error for $\psi_{6,4}$ decays down to $8\cdot 10^{-4}$ with decreasing mesh width, but then starts to rise again for some unknown reason.

\subsection{Cylinder flow}

As a second example, we reconsider a genuinely nonautonomous system with nonempty boundary\cite{FrLlSa10a}: a flow on the cylinder $M=S^1\times [0,\pi]$, given by the vector field
\begin{eqnarray*}
\dot x &=& c-A(t)\sin(x-\nu t)\cos(y) + \eps G(g(x,y,t))\sin(t/2)\\
\dot y &=& A(t)\cos(x-\nu t)\sin(y)
\end{eqnarray*}
with $A(t)=1+0.125\sin(2\sqrt{5}t)$, $G(\psi)=1/(\psi^2+1)^2$, $g(x,y,t)=\sin(x-\nu t)\sin(y)+y/2-\pi/4$ and parameter values $c = 0.5, \nu = 0.25$ and $\eps = 0.25$. We consider the flow map $T:=\Phi(\cdot,t_0,t_\mathbf{f})$ with $t_0=0$ and $t_\mathbf{f}=40$ which we approximate by 400 steps of the classical Runge-Kutta scheme of 4th order.

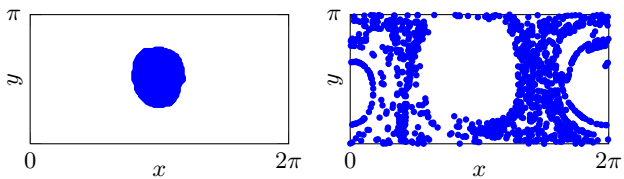
\begin{figure}[htbp]
  \centering
  \begin{tikzpicture}
    \begin{axis}[enlargelimits=false,
    	width=0.28\textwidth,
    	ylabel={$y$},ylabel style={at={(0.3,0.5)}},
    	xlabel={$x$},xlabel style={at={(0.5,0.2)}},
		xtick={0,6.28},ytick={0,3.14},
		xticklabels={0,$2\pi$},yticklabels={\empty,$\pi$},
		xmin=0,xmax=6.28,
		ymin=0,ymax=3.14,
		axis equal image,	
		mark size = 1
    ]
    \addplot[blue,only marks] table {cyl_non_coh_set.tex};
    \end{axis}
  \end{tikzpicture}
  \begin{tikzpicture}
    \begin{axis}[enlargelimits=false,
    	width=0.28\textwidth,
    	ylabel={$y$},ylabel style={at={(0.3,0.5)}},
    	xlabel={$x$},xlabel style={at={(0.5,0.2)}},
		xtick={0,6.28},ytick={0,3.14},
		xmin=0,xmax=6.28,
		ymin=0,ymax=3.14,
		xticklabels={0,$2\pi$},yticklabels={\empty,$\pi$},
		axis equal image,	
		mark size = 1
    ]
    \addplot[blue,only marks] table {cyl_non_coh_set_img.tex};
    \end{axis}
  \end{tikzpicture}

  \caption{Cylinder flow: Non-coherent pair. The set to the left is the one from Fig.~\ref{fig:cyl_curves} (left), shifted by $+2$ in $x$-direction.}\label{fig:cyl_non_coh_set}
\end{figure}

\begin{figure}[H]
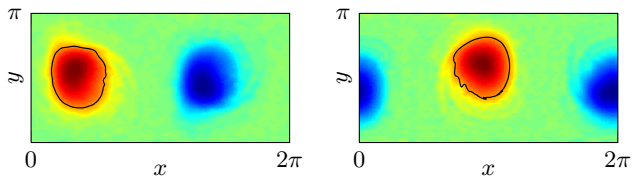

  \centering
  \cylevplot{cyl_ev2_curve}{}{$x$}{$y$}
  \cylevplot{cyl_ev2_curve_img}{}{$x$}{$y$}
  \caption{Minimising curves before (left) and after (right) the application of $\Phi(\cdot,t_0,t_\mathbf{f})$.}\label{fig:cyl_curves}
\end{figure}

The flow exhibits strong mixing apart from two eddy-like structures which roughly retain a constant $y$ value while moving along the circle (i.e.\ the $x$-) direction.  Figure~\ref{fig:cyl_curves} (left) shows one of these two coherent sets as a black curve, together with its image under the flow $\Phi(\cdot,t_0,t_\mathbf{f})$ (right).
On the other hand, selecting a set away from the objects identified in Fig.~\ref{fig:cyl_curves}, such as in Fig.~\ref{fig:cyl_non_coh_set}, reveals chaotic motion.  \OJ{There, the set within the black curve from Fig.~\ref{fig:cyl_curves} is shifted by 2 in $x$-direction at $t_0=0$, and shown (Fig.\ \ref{fig:cyl_non_coh_set}, left) together with its image at $t_\mathbf{f}=40$ (right).  Evidently, this  pair is not coherent.}

We applied Algorithm 1 with the following parameters: 1. For $\delta=10^{-6}$, we choose $Y=Y^\delta = ((h^\delta_x\Z+\delta)\times h_y\Z)\cap (S^1\times [0,\pi])$ with $h^\delta_x=(2\pi-2\delta)/50$ and $h_y=\pi/50$ for the centers, i.e.\ $2500$ points on a regular grid which is slightly shifted away from the boundary, and $\eps=2$ for the shape parameter such that the supports of the basis functions overlap with roughly 100 others.  2. We choose $X_\inn  = Y^0\backslash \{(x,y)\mid y=0\text{ or } y = \pi\}$ for the inner collocation points and $X_\bd  = Y^0\backslash \{(x,y)\mid y\neq 0\text{ and } y \neq \pi\}$ for the boundary points.  The Cauchy-Green tensor was computed by integrating the (backward) variational equation for each point in $X_\bd $.  3. We compute the matrices $P$ and $D$ as described in Sections~\ref{subsec:first_case} and \ref{sect:implementation}
and 4. solve the eigenproblem (\ref{eq:discrete_eigenproblem}).

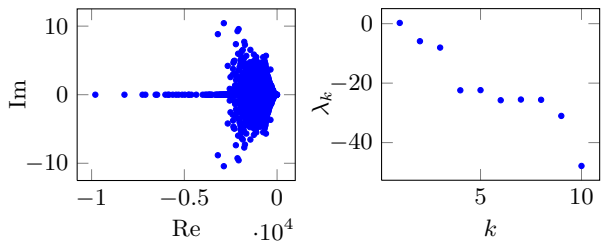
\begin{figure}[htbp]
  \centering
   \begin{tikzpicture}
    \begin{axis}[
    	width=0.25\textwidth,
    	ylabel={Im},ylabel style={at={(0.15,0.5)}},
    	xlabel={Re},xlabel style={at={(0.5,0.05)}},
		mark size = 1
    ]
    \addplot[blue,only marks] table {cyl_specDL.tex};
    \end{axis}
  \end{tikzpicture}
   \begin{tikzpicture}
    \begin{axis}[
    	width=0.25\textwidth,
    	ylabel={$\lambda_k$},ylabel style={at={(0.15,0.45)}},
    	xlabel={$k$},xlabel style={at={(0.5,0.05)}},
		mark size = 1
    ]
    \addplot[blue,only marks] table {cyl_specDL_close.tex};
    \end{axis}
  \end{tikzpicture}
  \caption{Cylinder flow: Spectrum of $\frac12 (D+P^\top DP)$ (left: entire spectrum, right: eigenvalues closest to~$0$).}\label{fig:cyl_spectrum}
\end{figure}

The resulting four eigenvalues with smallest magnitude were $\lambda_1=0.25$, $\lambda_2 = -5.90$, $\lambda_3=-8.10$ and $\lambda_4=-22.4$.
Fig.~\ref{fig:cyl_spectrum} shows the spectrum of the discrete dynamic Laplacian $\frac12 (D+P^\top DP)$.  Again, the spectrum is not real, but it is real close to $0$ which is the part we are interested in (cf.\ our comments on the spectrum in Section \ref{sec:standard}).

\begin{figure}[htp]
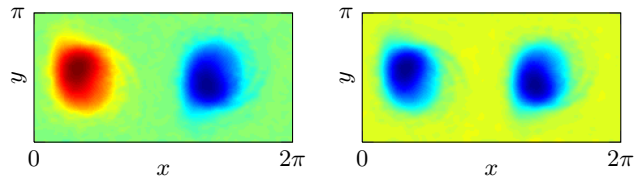

  \centering
  \cylevplot{cyl_ev1}{}{$x$}{$y$}
  \cylevplot{cyl_ev2}{}{$x$}{$y$}
  \caption{Cylinder flow map: Eigenvectors $f_2$ (left) and $f_3$ (right) of (\ref{eq:discrete_eigenproblem}), using $n=2500$ centers on a regular grid with $\eps=2$. }\label{fig:cyl_evecs}
\end{figure}

\begin{figure}[htp]
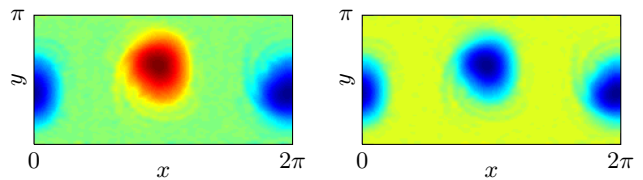

  \centering
  \cylevplot{cyl_ev1_img}{}{$x$}{$y$}
  \cylevplot{cyl_ev2_img}{}{$x$}{$y$}
  \caption{Images $Pf_2$ and $Pf_3$ of the eigenvectors from Figure \ref{fig:cyl_evecs}.}\label{fig:cyl_pushevecs}
\end{figure}

We again see that $\lambda_2$ and $\lambda_3$ are relatively close, compared to lower eigenvalues.
This indicates the presence of two distinct finite-time coherent sets.
The corresponding eigenfunctions $f_2$ and $f_3$ are shown in Figure~\ref{fig:cyl_evecs}, while Figure~\ref{fig:cyl_pushevecs} shows their push-forwards $Pf_2$ and $Pf_3$.

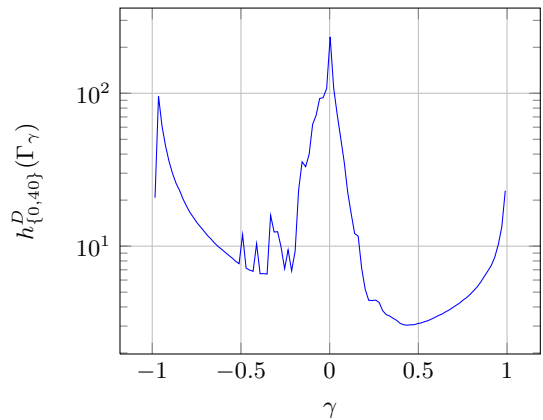
\begin{figure}[htbp]
  \centering
   \begin{tikzpicture}
    \begin{semilogyaxis}[
    	width=0.4\textwidth,
    	ylabel={$h^D_{\{0,40\}}(\Gamma_\gamma)$},
    	xlabel={$\gamma$},
    	grid
    ]
    \addplot[blue] table {cyl_cheeger.tex};
    \end{semilogyaxis}
  \end{tikzpicture}

  \caption{Cylinder flow: Graph of the dynamic Cheeger constant (\ref{cheegereqn}) in dependence of the level set value $\gamma$.}\label{fig:cyl_thresh}
\end{figure}

Again, as described in Step 6 of Algorithm 1, we scan over the level sets of $f_2$ (shown to the left in Figure \ref{fig:cyl_evecs}) and compute the length of the level set and the length of its image.
To illustrate the first of the two options in Step 6 of Algorithm 1, we now use Step 6(a) and directly map points found on the curves $\Gamma_\gamma=\{x\in S^1\times [0,\pi]: f_2(x)=\gamma\}$,  $\gamma\in[-1,1]$, to evaluate $\mathbf{h}^D_{\{0,40\}}(\Gamma_\gamma)$.
Figure \ref{fig:cyl_thresh} shows the value of (\ref{cheegereqn}) as a function of $\gamma$.
The minimising $\gamma=0.4372$ was chosen by scanning through 100 equally separated values of $\gamma$ from $\min f_2$ to $\max f_2$, where $f_2$ is scaled so that $|f_2|_\infty=1$.
The curve along the corresponding level set in $f_2$ and its image under $\Phi(\cdot,0,40)$ are shown in Figure \ref{fig:cyl_curves}, having lengths $\ell_1(\Gamma_\gamma)=4.52$ and $\ell_1(\Phi(\Gamma_\gamma,t_0,t_\mathbf{f}))=5.01$, while $\ell_2(M_1)=1.57$ for the corresponding finite-time coherent set $M_1$ (red in the Figure) so that $h^D_{0,40}(\Gamma_\gamma)=3.04$.
One could now threshold $f_3$ according to Steps 6 and 7 of Algorithm 1 to find an optimal $\Gamma$ from $f_3$.
This would yield the union of the curve already found from $f_2$ and a closed curve surrounding the right-hand blue set in Figure \ref{fig:cyl_curves}.

\section{Conclusion and future directions}

\vspace*{-2mm}

We presented a fast method for the discretization of transfer operators based on collocation with radial basis functions (RBFs), and applied this method to compute finite-time coherent sets using the new advection-only construction from \cite{F15}.
In particular, by choosing sufficiently smooth kernels, we observed that very high convergence orders can be achieved.
These rapid constructions of accurate numerical approximations of transfer operators alleviate the major computational expense in algorithms to detect finite-time coherent sets.

We demonstrated that we could construct accurate estimates of the boundaries of finite-time coherent sets in a nonlinear cylinder flow over a rather long time duration using only a $50\times 50$ grid of initial and final points, and without exploiting any dynamics or geometry of the flow in our choice of RBF centres or radii.
Beyond constructing a discretisation of the main boundary value problem, our numerical contributions include (i) stable methods to handle the boundary conditions that avoid matrix inversion, and (ii) methods to deal with the situation where the entire domain evolves under the dynamics.

In our numerical experiments we observed that the accuracy of the RBF collocation procedure can be sensitive to the various parameters involved, i.e.\ to the choice of the centers, the collocation points, and the shape parameter $\eps$, which controls the radius of the supports of the basis functions.  In practice, it is therefore advisable to repeat a particular computation with different values for the parameters.  We observed that choosing the collocation points different from the centers does not yield a real spectrum of the discrete Laplacian. A more systematic study of this phenomenon will be the subject of future work.

\bibliographystyle{abbrv}

\end{document}